\documentclass[a4paper,12pt]{article}
\usepackage{latexsym}
\usepackage{amsmath}
\usepackage{amsfonts}
\usepackage{amssymb}
\usepackage{graphicx}
\usepackage{color}
\usepackage[LGR,T1]{fontenc}

\begin{document}

\title{Exact soliton-like probability measures for interacting jump processes}

 \author{ $\,$\\M.-O. Hongler, (max.hongler {\it at} epfl.ch)\\ EPFL/ IMT/ LPM \\ Station 17 \\ CH-1015 Lausanne, (Switzerland) }
\date{}
\maketitle

\abstract{

\noindent  The cooperative dynamics  of a 1-D collection of  Markov jump, interacting stochastic processes is studied via a mean-field approach. In the time-asymptotic regime, the resulting  nonlinear master equation is analytically solved. The nonlinearity compensates  jumps induced diffusive behavior giving rise to a soliton-like stationary probability density. The soliton velocity and its sharpness both intimately depend on the interaction strength. Below a critical threshold of the strength of interactions, the cooperative behavior cannot be sustained leading to the destruction of the  soliton-like solution.  The bifurcation point for this behavioral phase transition  is explicitly calculated.}

\vspace{0.2cm}

\vspace{0.2cm}
\noindent {\bf Classification}: 60J75

\section{Introduction}

\noindent Interacting stochastic agents are modeled by a collection of  nonlinearly coupled Markovian   stochastic processes. Inspired by  the dynamics recently exposed in [Bal\'azs 2014],  we  focus  on pure, right-oriented  jump processes.  For large and homogeneous swarms, the mean-field  description offers a powerful method to characterize  the resulting nonlinear global dynamics.  Adopting the MF approach, the swarm behavior is summarized into a field density variable obeying a nonlinear master equation.  Such  partial differential  integral equations  are in general   barely   completely solvable. Nevertheless, several explicitly solvable models have been recently studied [Hongler 2014, Bal\'azs 2014]. Our present goal is  to enrich this yet available  collection by  proposing an intrinsically nonlinear extension of the recent models introduced by  [Bal\'azs 2014]. Models involving pure jumps  complete  the solvable models with dynamics driven either by Brownian Motion  and or by  alternating Markov renewal processes [Hongler 2014]. For strong enough mutual  interactions, we explicitly observe the existence of a stationary probability measure  propagating like a soliton. This  soliton-like  dynamics can be formed  since the underlying nonlinear mechanism due to interactions exactly  compensates the  jump induced diffusion.  This exhibits a  close analogy with  nonlinear wave dynamics where nonlinearity compensates the  velocity dispersion.  Since the model is uni-dimensional, long-range interactions between the agents are mandatory for the existence of cooperative behaviors here described by  soliton-like probability measures.  Decreasing the strength of the mutual  interactions, via a barycentric modulation function similar to the one used in [Bal\'azs 2014], we reach a critical threshold below which no stable cooperative behavior can be sustained. The critical threshold where the  behavioral  phase transition occurs can here be exactly calculated.

\section{Linear pure jump stochastic processes}

\noindent Let us first describe the dynamics of a single, isolated jump  process which later in section 2, will  enter into the composition of our interacting swarm. On $\mathbb{R}$, we consider the  right-oriented jump  Markovian process  $X(t)$ characterized by  the (linear) master equation: 

\begin{equation}
\label{BASE/LIN}
\partial_{t}P(x,t) =  - P(x,t) + \int_{-\infty}^{x}  P(y,t)  \varphi( x-y) dy,
\end{equation}

\noindent where $P(x,t)$ with $P(x,0) = f(x)$ stands for the transition probability density. The function $\varphi(x) : \mathbb{R} \rightarrow \mathbb{R}^{+} $ defines the probability density for the (right oriented) lengths of the process jumps.

\noindent Taking the $x$-Laplace transform of Eq.(\ref{BASE/LIN}) and taking into account of the convolution structure, we obtain directly:

\begin{equation}
\label{CONVOLIN}
\partial_t \tilde{P}(s, t) = - \left[ 1 - \tilde{\varphi}(s) \right] \tilde{P}(s,t).
\end{equation}

\noindent which solution reads:

\begin{equation}
\label{EVOLIN}
\tilde{P}(s, t) = e^{-  t +  \tilde{\varphi}(s) t },
\end{equation}

\noindent where in writing Eq.(\ref{EVOLIN}), we already did  assume the initial condition:

\begin{equation}
\label{INIT}
P(x,t)\mid _{t=0}= \delta(x).
\end{equation}

\vspace{1cm}
\noindent {\bf Example}. Consider the dynamics obtained when $\varphi(x) = \lambda e^{- \lambda x}$ yielding  $\tilde{\varphi}(s) = {\lambda \over \lambda +s}$ and when the initial probability density is $f(x) = \delta(x)$.
Accordingly Eq.(\ref{EVOLIN}) reads:
\begin{equation}
\label{EVOLINK}
\tilde{P}(s,t) = e^{-  t} \left[e^{ t \left( {\lambda \over \lambda +s }\right)} \right] =  e^{-  t} \left\{\sum_{n=0}^{\infty}  {(\lambda \,   t )^{n} \over n!} \left[ {1 \over \lambda +s }\right]^{n}\right\}.
\end{equation}

\noindent The Laplace inversion of Eq.(\ref{EVOLINK}) yields:

\begin{equation}
\label{EVOLM}
P(x,t) =  e^{-  t} \left\{\delta(x)  +  e^{- \lambda x}\underbrace{ \sum_{n=1 }^{\infty}  {(\lambda \,   t )^{n} \over n!}  \, { x^{n-1}\over (n-1)!}}_{:= J(x,t)}
\right\}.
\end{equation}

\noindent For  $J(x,t)$, we can write:

\begin{equation}
\label{SUMJ}
J(x,t) = {d \over dx } \left\{ \underbrace{ \sum_{n=1 }^{\infty}{(\lambda   t )^{n} \over n!}  \, { x^{n} \over n !}}_{\mathbb{I}_{0} \left(2 \sqrt{\lambda  \, x\, t} \right)-1}\right\}  =  {\sqrt{\lambda t} \over \sqrt{x} } \mathbb{I}_1\left(2 \sqrt{\lambda  \, x\, t} \right)\end{equation}

\noindent where  $\mathbb{I}_{m}(z)$ stands for the $m$-modified Bessel's functions. 

\noindent Hence  the final probability density $ P(x,t$) reads:

\begin{equation}
\label{QPRO}
P(x,t) =  e^{-  t } \left\{ \delta(x)  + e^{- \lambda x}{\sqrt{\lambda  t} \over \sqrt{x} } \mathbb{I}_1\left(2 \sqrt{\lambda \,x\, t} \right)
\right\}, \quad x \in \mathbb{R}^{+}
\end{equation}

\noindent  and one may explicitly verify that one indeed has;  $\int_{\mathbb{R}^{+}} P(x,t)  dx=1$, (use the  entry  6.643(2) in [Gradshteyn 80]).

\noindent For time asymptotic regimes, Eq.(\ref{QPRO}) behaves as:

\begin{equation}
\label{QAS}
\lim_{t \rightarrow \infty} P(x,t) \simeq {(\lambda\,   t)^{{1 \over 4} }\over 2 \sqrt{\pi} x^{{3 \over 4}}} e^{ -\left[\sqrt{\lambda x}  - \sqrt{  t }\right]^{2}},
\end{equation}

\noindent  exhibiting therefore a  diffusive  propagating wave  with vanishing amplitude and  velocity $V:= {1  \over \lambda}$. Due to translation invariance of the dynamics, we note that  $P(x- y,t)$ fulfills a $\delta(x-y)$ initial condition. Hence,  when $P(x,0) = f(x)$, the linearity of the dynamics Eq.(\ref{BASE/LIN})  enables to write:

\begin{equation}
\label{ARBOS}
\left\{ 
\begin{array}{l}
P_f(x,0) = f(x), \\  \\ 
P_f(x,t) = \int_{\mathbb{R}^{+}} P( (x-y) ,t) f(y) dy. 
\end{array}
\right.
\end{equation} 

\section{Non-linear Markovian  jump processes}

\noindent Keeping the jumps probability density as $\varphi(x) = \lambda e^{- \lambda x}$, let us now consider  a large homogeneous collection of identical  processes evolving like Eq.(\ref{BASE/LIN}) now subject of  mutual long-range  interactions. The class of interactions we consider  yields,  in the  mean-field limit,  the nonlinear master equation:

\begin{equation}
\label{KOMP1}
\left\{ 
\begin{array}{l}
\Omega(x,t) =   \int_{x}^{\infty}   g\left(z - \langle X(t) \rangle \right) \partial_{z}G(z,t) dz\\ \, \\ 
\partial_{xt}G(x,t) =  -\Omega(x,t) \partial_{x}G(x,t)  + \int_{-\infty}^{x}  \Omega (y, t)  \partial_{y}G(y,t)  \lambda e^{-\lambda (x-y)} dy, \\ \, \\
\langle X(t) \rangle = \int_{\mathbb{R}^{+} }y \, \partial_{y}G(y,t) \, dy,
\end{array}
\right.
\end{equation}

\noindent where $G(x,t)$ stands for the cumulative distribution of the a nonlinear jump process, (i.e. $G(x,t)$ is monotonously increasing with boundary conditions  $G(-\infty,t)=0$ and $G(\infty,t)=1$). Note that while in  Eq.(\ref{BASE/LIN})  the  jumping rate is unity,  in Eq.(\ref{KOMP1}) it is replaced by   $\Omega(x,t)>0$ which is explicitly state-dependent. This is precisely where  the mutual interaction introduce  a strong nonlinearity  into the dynamics. In the sequel, we  focus on cases where  $g(x) = g(-x) >0$.

\vspace{0.2cm}
\noindent For asymptotic time, we now postulate that Eq.(\ref{KOMP1}) admits $\xi$-functional dependent solutions with  $\xi = (x - Vt)$ and  with the even symmetry:

\begin{equation}
\label{AVERAGE}
\int_{\mathbb{R}} \xi \partial_{\xi} G(\xi)  d\xi  =0, 
\end{equation}

\noindent where $V$ is a propagating velocity parameter. In terms of $\xi$, Eq.(\ref{KOMP1}) can be rewritten as:

\begin{equation}
\label{TRAVELO1}
 V \left[ \partial^{3}_{\xi\xi \xi} G(\xi) +  \lambda \partial^{2}_{\xi\xi } G(\xi) \right] = \partial_{\xi} \left\{\Omega (\xi)  \partial_{\xi} G(\xi) \right\}.
\end{equation}

\noindent Defining  ${\cal L} (\xi) := \log \left[\partial_{\xi} G(\xi)  \right]$, after one integration step where the integration constant is taken to  be zero, Eq.(\ref{TRAVELO1}) can be rewritten as:

\begin{equation}
\label{TRAVELO2}
V \partial^{2}_{\xi\xi} {\cal L} (\xi)  =  -\lambda  V +  \int_{\xi}^{\infty} \left[ g(\eta) \partial_{\eta} G(\eta )\right] d \eta.
\end{equation}

\noindent Assuming now a  functional dependence $g(\xi) =  \cosh^{-n}( \xi)$ with    $n \in \mathbb{R}$, by direct substitution,  it is immediate to see that  Eq.(\ref{TRAVELO2}) is  solved by  the  (normalized) probability density $\partial_{\xi} G(\xi) $:

\begin{equation}
\label{ANSATZOS}
\left\{
\begin{array}{l}
\partial_{\xi} G(\xi) = { \Gamma \left( {m+1 \over 2} \right) \over \sqrt{\pi} \Gamma \left( {m \over 2} \right) } \cosh ^{-m} ( \xi),  \qquad  
m>0, \\ \, \\ m=  \lambda = 2-n, \\ \, \\
 V =  {\Gamma\left( {m+1\over 2} \right) \over  \sqrt{\pi} \Gamma \left( {m \over 2} \right) } .
\end{array}
\right.
\end{equation}

\noindent Due to the $\xi$-symmetry of the probability density $\partial_{\xi} G(\xi) $, Eq.(\ref{AVERAGE}) is trivially fulfilled.

\noindent For $n \in ]2,-\infty]$, Eq.(\ref{ANSATZOS}) implies  that a stationary propagating  density $\partial_{x}G(x)$ is sustained by the nonlinear dynamics Eq.(\ref{KOMP1}). However, for short decaying $g(x)$-modulation, occurring when  $n>2$, no stationary propagating probability density exists, (i.e. for this parameter range,  $m<0$ in Eq.(\ref{ANSATZOS}) and the solution  cannot be normalized to unity as required for a probability measure).  For this exactly solvable dynamics, we also observe that the average jump length $\lambda^{-1}$ and the barycentric modulation strength controlled by the factor  $n$ are intimately dependent control  parameters. In addition,  we note that for large $m$, the asymptotic expansion of the $\Gamma$-function implies that $\lim_{m \rightarrow \infty}V \simeq \sqrt{m}$.

\vspace{0.4cm}
\noindent {\bf Illustration}. Along the same lines as exposed in  [Hongler 2014],  the nonlinear dynamics given by Eq.(\ref{KOMP1})  can be viewed as representing  the mean-field evolution associated with a large population of  stochastic  jumping agents  subject to a mutual imitation process. The swarm dynamics is described via the  probability density function $\partial_xG(x,t)$ obeying a nonlinear partial differential equations (PDE).  Agents mutual interactions are responsible for the state-dependent jumping rate $\Omega(x,t)$ in  Eq.(\ref{KOMP1}). The functional form of $\Omega(x,t)$   simultaneously englobes two distinct nonlinear features, namely:

\begin{itemize}
  \item[]
  \item[] a) {\bf imitation process}.  To isolate this process, we may consider the case $g(x) \equiv 1$, (i.e. $n=0$) implying that   
  
  \begin{equation}
\label{IMITATION}
\Omega(x,t) = 1 - G(x,t).
\end{equation}
The resulting  state-dependent jumping rate Eq.(\ref{IMITATION}) induces a traveling and compacting tendency. As the agents are subject to pure right-oriented jump, Eq.(\ref{IMITATION}) effectively describes situations where the laggard agents  jump more frequently than  the leaders, (i.e. laggards try to effectively imitate the leaders behavior).
\item[]
  \item[] b) {\bf barycentric range modulation of the mutual  interactions}. The modulation obtained  when $g(x) \neq 1$ describes  the relative importance attributed to interactions with agents  remote  from  the barycenter $ \langle X(t) \rangle $ of the  swarm. Here, we may separate two distinct  tendencies: 
  
  \begin{itemize}
  \item[]
  \item[] i)  when $n \in [0, 2[$, far remote agents tend to  not influence the dynamics. In this case, the resulting behavior can be referred as a  {\bf weak cooperate identity} and the propagating probability density  given by Eq.(\ref{ANSATZOS}) exhibits the shape of {\bf a table-top soliton} with  a plateau increasing when the limiting value $2$ is approached. AOn observes a comparatively  low propagating velocity $V$ of these table-top like aggregates. Again, we emphasize  that for $n>2$, the cooperative interactions are not strong enough to sustain the propagation  of a cooperative behavior in asymptotic  time. This is well known in general  for 1-D stochastic interacting system, (the Ising model being the paradigmatic example) where no cooperative phase can be formed  when the interactions operate on too limited ranges.
  \item[]
  \item[] ii) for $n < 0$, the $g(x)$ modulation effectively gives rise to a {\bf strong cooperate identity}. Far remote agents increasingly  influence  the swarm. This gives rise to sharply peaked solitons-like probability densities propagating with high propagating velocities.
\end{itemize}
  
  \item[] 
\end{itemize}
\noindent

\section*{References}

\noindent  [Bal\'azs 2014]. M. Bal\'azs, M. Z. Mikl\'os and B. T\'oth. {\it "Modeling flocks and prices: Jumping particles with an interactive interaction"}. Ann. Inst. H. Poincar\'e, {\bf50}(2), (2014), 425-454.

\vspace{0.2cm}
\noindent [Hongler 2014]. M.-O. Hongler, R. Filliger and O. Gallay. {\it "Locas versus nonlocal barycentric interactions in 1D agent dynamics"}. Mathematical Biosciences {\bf 11}(2), April 2014.

\vspace{0.2cm}
\noindent [Gradshteyn 80]. I. S. Gradshteyn and I. M. Ryzhik. {\it Table of Integrals, Series and Products}. Academic press (1980).

\end{document}